\documentclass[12pt,a4paper]{article}
\usepackage[utf8]{inputenc}
\usepackage{amsmath}
\usepackage{amsfonts}
\usepackage{amssymb}
\usepackage{amsthm}
\usepackage{authblk}
\usepackage{inputenc}
\usepackage{longtable}

\allowdisplaybreaks

\begin{document}
\title{A New Type of Abundant Numbers}
\author{Xiaolong Wu}
\affil{Ex. Institute of Mathematics, Chinese Academy of Sciences}
\affil{xwu622@comcast.net}
\date{Jun 13, 2019}
\maketitle

\begin{abstract}

    This article defines a new type of abundant numbers, called largest rho-value (abbreviate LR) numbers, and then shows that Robin hypothesis is true if and only if all LR numbers $>5040$ satisfy Robin inequality. 
\end{abstract}
\begin{center}\textbf{ \large Introduction}
\end{center}

Let $n\geq 2$ be an integer, $\sigma (n)=\sum_{d|n}$ is the sum of divisor function.  Define
\begin{equation}
\rho (n):=\frac{\sigma (n)}{n},
\end{equation}
\begin{equation}
G(n):=\frac{\rho (n)}{\log \log n},
\end{equation}

    Robin [Robin 1984] made hypothesis that all integers $n>5040$ satisfy Robin inequality
\begin{equation}\tag{RI}
G(n)<e^\gamma ,
\end{equation}
where $\gamma$ is the Euler constant.

Write the factorization of n as
\begin{equation}
n=\prod_{i=1}^r p_i^{a_i},
\end{equation}
where $p_i$ is the i-th prime and $a_r\geq 1$. If $p_i$ and n are co-prime, we set its exponent $a_i=0$.

    Define the function of sum of exponents:
\begin{equation}
\omega (n) :=\sum_{i=1}^r a_i.
\end{equation}

For an integer $m\geq 1$, define a set
\begin{equation}
S_m :=\{n\in \mathbb{N} \, |\, n\geq 2, \,\omega (n)=m \}.
\end{equation}

    If there exists an element $n_m\in S_m$ such that $\rho(n_m )\geq \rho(n),\, \forall \,n \in S_m$ , we call $n_m$ a largest rho-value (abbreviate LR) number. 
    
 We now construct $n_m$. Define a set
\begin{equation}\notag
Z :=\{z=z_{q,k}:=q+\cdots +q^k \, | \, q \, prime,\quad k\geq1 \, integer \}.
\end{equation}
We sort Z in increasing order. 
    
        Note that Z contains duplicate integers. One known example is $5+5^2=2+2^2+2^3+2^4$. If $q+\cdots +q^k=p+\cdots +p^j$ and $q>p$, then we assign $z_{q,k}$ smaller order than $z_{p,j}$. We will denote by $z_i$, $q_i$ and $k_i$ the numbers related to the i-th element. 

Next, we define
\begin{equation}
\delta_i=\delta_{q_i,k_i} :=\log \left(1+\frac{1}{z_i}\right)=\log \left(1+\frac{1}{q_i+\cdots +q_i^{k_i}}\right).
\end{equation}
Since $z_i$ are increasing, $\delta_i$ are decreasing. Define the set of triplets
\begin{equation}
E:=\{ (q,k,\delta_{q,k})\}.
\end{equation}
Elements in E are ordered by $\delta_i$.

    Now for a given integer $m\geq 1$, let $E_m$ be the subset of first m elements in E. For any prime p, let $k_p$ be the largest k among elements $\delta_{p,k}\in E_m$. Define
\begin{equation}
n_m:=\prod_p p^{k_p}.
\end{equation}
    
    Theorem 1 will show that $n_m$ is an LR number. That is, $\rho(n_m )\geq \rho(n),\, \forall \,n \in S_m$.
    
\begin{center}
             \textbf{Table 1. First 20 LR numbers}
\begin{longtable}{| r | r | r | r | r | r | r | r | r |}
\hline
\textbf{m} & \textbf{q} & \textbf{k} & $z_m$ & $\delta_{q,k}$ & $n_m$ &$\rho(n_m)$&$\log(n_m)$&$G(n_m)$\\
\hline
1	&2	&1	&2	&0.4055	&2	&1.5000	&0.6931	&-4.0296\\
\hline
2	&3	&1	&3	&0.2877	&6	&2.0000	&1.7918	&3.4294\\
\hline
3	&5	&1	&5	&0.1823	&30	&2.4000	&3.4012	&1.9606\\
\hline
4	&2	&2	&6	&0.1542	&60	&2.8000	&4.0943	&1.9864\\
\hline
5	&7	&1	&7	&0.1335	&420	&3.2000	&6.0403	&1.7793\\
\hline
6	&11	&1	&11	&0.0870	&4,620	&3.4909	&8.4381	&1.6368\\
\hline
7	&3	&2	&12	&0.0800	&13,860	&3.7818	&9.5368	&1.6770\\
\hline
8	&13	&1	&13	&0.0741	&180,180	&4.0727	&12.1017	&1.6334\\
\hline
9	&2	&3	&14	&0.0690	&360,360	&4.3636	&12.7949	&1.7119\\
\hline
10	&17	&1	&17	&0.0572	&(17\#)(3\#)2	&4.6203	&15.6281	&1.6807\\
\hline
11	&19	&1	&19	&0.0513	&(19\#)(3\#)2	&4.8635	&18.5725	&1.6646\\
\hline
12	&23	&1	&23	&0.0426	&(23\#)(3\#)2	&5.0750	&21.7080	&1.6490\\
\hline
13	&29	&1	&29	&0.0339	&(29\#)(3\#)2	&5.2500	&25.0753	&1.6295\\
\hline
14	&5	&2	&30	&0.0328	&$(29\#)(5\#)2$	&5.4249	&26.6847	&1.6519\\
\hline
15	&2	&4	&30	&0.0328	&$(29\#)(5\#)2^2$	&5.6058	&27.3779	&1.6937\\
\hline
16	&31	&1	&31	&0.0317	&$(31\#)(5\#)2^2$	&5.7866	&30.8119	&1.6881\\
\hline
17	&37	&1	&37	&0.0267	&$(37\#)(5\#)2^2$	&5.9430	&34.4228	&1.6794\\
\hline
18	&3	&3	&39	&0.0253	&(37\#)(5\#)(3\#)2	&6.0954	&35.5214	&1.7073\\
\hline
19	&41	&1	&41	&0.0241	&(41\#)(5\#)(3\#)2	&6.2441	&39.2350	&1.7016\\
\hline
20	&43	&1	&43	&0.0230	&(43\#)(5\#)(3\#)2	&6.3893	&42.9962	&1.6988\\
\hline
\end{longtable}
\end{center}
Here “\#” means primorial.

    Theorems 2-8 study properties of $n_m$. Theorem 9 shows that Robin hypothesis is true if and only if all LR numbers $n_m>5040$ satisfy Robin inequality.
    
    We will use $\theta(x)$and $\psi(x)$ for Chebyshev functions.

    For an integer $m\geq 1$ define $y_{m,1}:=z_m$, and for $k\geq 2$, define $y_{m,k}$ as the solution of
\begin{equation}
y_{m,k}+\cdots +y_{m,k}^k=z_m.
\end{equation}
When m is obvious, we will simply write $y_k$ instead of  $y_{m,k}$.
    
\begin{center}
\textbf{ \large Main Content}
\end{center}

\noindent {\bfseries Theprem 1.}
\textit{ Let $n_m$ be constructed as in (8). Then $\rho(n_m )\geq \rho(n),\, \forall \,n \in S_m$. That is, $n_m$ is an LR number.}
\begin{proof}
For a prime p, we have $\log \rho(p)=\log (1+1/p)=\delta_{p,1}$, and for integer $k\geq 2$, 
\begin{align*}
\log \rho(p^k)-\log \rho(p^{k-1})&=\log \frac{1+\cdots +p^k}{p^k}-\log \frac{1+\cdots +p^{k-1}}{p^{k-1}}\\
&=\log \left(\frac{1+\cdots +p^k}{p(1+\cdots +p^{k-1})}\right)\\
&=\log \left(1+\frac{1}{p+\cdots +p^k}\right)=\delta_{p,k}.\tag{1.1}
\end{align*}
Hence
\begin{equation}\tag{1.2}
\log \rho(p^k)=\sum_{i=1}^k \delta_{p,k},\quad \forall \, k\geq1.
\end{equation}
For any integer
\begin{equation}\notag
n=\prod_p p^{a_p}
\end{equation}
in $S_m$ we have $\sum_p a_p=m$ and
\begin{equation}\tag{1.3}
\log \rho(n)=\log \prod_p \rho(p^{a_p})=\sum_p\log \rho(p^{a_p})=\sum_p \sum_{i=1}^{a_p} \delta_{p,i}.
\end{equation}
By construction of $n_m$, $\log n_m=\sum_p \sum_{i=1}^{k_p} \delta_{p,i}$ is the largest sum of m elements in E. So $\log \rho(n_m )\geq \log \rho(n)$ for all $n\in S_m$. Hence $\rho(n_m )\geq \rho(n)$ or all $n\in S_m$. 
\end{proof}

\noindent {\bfseries Lemma 1.}
\textit{For given $z_m$ and $k\geq 2$, we have 
\begin{equation}\tag{L1.1}
z_m^{1/k}-1<y_k<z_m^{1/k}.
\end{equation}
\begin{equation}\tag{L1.2}
\theta(z_m^{1/k})-\frac{1}{k}\log z_m\leq \theta(y_k)\leq \theta(z_m^{1/k}).
\end{equation}
}
\begin{proof}
\begin{equation}\tag{L1.3}
(y_k+1)^k>y_k+\cdots y_k^k=z_m\,\Longrightarrow \, y_k>z_m^{1/k}-1.
\end{equation}
\begin{equation}\tag{L1.4}
y_k^k<y_k+\cdots y_k^k=z_m\,\Longrightarrow \, y_k<z_m^{1/k}.
\end{equation}
Since there is at most one prime in $(z_m^{1/k}-1,\,z_m^{1/k}]$, we have
\begin{equation}\tag{L1.5}
\theta(y_k)>\theta(z_m^{1/k})-\log(z_m^{1/k})=\theta(z_m^{1/k})-\frac{1}{k}\log z_m.
\end{equation}
\end{proof}

\noindent {\bfseries Lemma 2.}
\textit{For given $n_m$, we have for each prime $p<z_m$
\begin{equation}\tag{L2.1}
\left\lfloor\frac{\log z_m}{\log p}\right\rfloor-1\leq k_p\leq \left\lfloor \frac{\log z_m}{\log p}\right\rfloor.
\end{equation}
}
\begin{proof}
From the construction of  $n_m$, $p^{k_p}\leq p+\cdots+p^{k_p}\leq z_m$ means $k_p\leq \log z_m /\log p$. \\
$p^{k_p+2}\geq 2p^{k_p+1}>(1+p^{-1}+\cdots+p^{-k_p})p^{k_p+1}=p+\cdots+p^{k_p+1}>z_m$ means $k_p+2>\log z_m/\log p$. Hence (L2.1) holds.
\end{proof}

\noindent {\bfseries Theorem 2.}
\textit{Let $m\geq 1$ be an integer, $n_m$ and $k_p$ be defined as in (8), $y_k$ be defined as in (9). Then
\begin{equation}\tag{2.1}
k_p=\begin{cases}
k, & y_{k+1}<p\leq y_k,\\
0, & y_1<p.
\end{cases}
\end{equation} 
when $\delta_m\neq \delta_{m+1}$.\\
When $\delta_{q,k}=\delta_m= \delta_{m+1}=\delta_{q',k'}$ and $q>q'$, then $k_q=k$, $k_{q'}=k'-1$.
}
\begin{proof}
Define
\begin{equation}\tag{2.2}
D(x,k):=\log \left(1+\frac{1}{x+\cdots +x^k}\right). 
\end{equation}
First assume $\delta_m\neq \delta_{m+1}$. Since $D(x,k)$ decreases in x, $p\leq y_k$ means $D(p,k)\geq D(y_k,k)=D(z_m,1)$. Hence $D(p,k)=\delta_{p,k}\in E_m$. That is, $k_p\geq k$. On the other hand, $y_{k+1}<p$ means $D(p,k+1)<D(y_{k+1},k+1)=D(z_m,1)$. Hence $D(p,k+1)=\delta_{p,k+1}\notin E_m$. That is, $k_p<k+1$. So $k_p=k$.

The $\delta_m=\delta_{m+1}$ case is obvious from the requirement "larger prime has smaller index" when sorting Z.
\end{proof}

Define two constants
\begin{equation}
W_1:=\sum_{i=1}^\infty \left(\frac{1}{z_i}-\log \left(1+\frac{1}{z_i}\right)\right)=0.20208\cdots,
\end{equation}
\begin{equation}
W_2:=M+\sum_{z\in Z,z\neq prime}\frac{1}{z}=0.77929\cdots,\quad\quad\quad
\end{equation}
here M is the Meissel-Mertens constant
\begin{equation}\notag
M:=\gamma+\sum_p \left(\log \left(1-\frac{1}{p}\right)+\frac{1}{p}\right)=0.26149\cdots.
\end{equation}

\noindent {\bfseries Theorem 3.}
\textit{
\begin{equation}\tag{3.1}
-W_1+W_2=\gamma.
\end{equation}
}
\begin{proof}
We have
\begin{align*}
-W_1&+W_2-\gamma\\
&=-\sum_{i=1}^\infty \left(\frac{1}{z_i}-\log \left(1+\frac{1}{z_i}\right)\right)\\
&\quad +\gamma +\sum_p \left(\log \left(1-\frac{1}{p}\right)+\frac{1}{p}\right)+\sum_{z\in Z,z\neq prime}\frac{1}{z}-\gamma\\
&=-\sum_p \left(\frac{1}{p}-\log \left(1+\frac{1}{p}\right)\right)-\sum_{z\in Z,z\neq prime} \left(\frac{1}{z}-\log \left(1+\frac{1}{z}\right)\right)\\
&\quad +\sum_p \left(\log \left(1-\frac{1}{p}\right)+\frac{1}{p}\right)+\sum_{z\in Z,z\neq prime}\frac{1}{z}\\
&=\sum_p \left(\log \left(1+\frac{1}{p}\right)+\log \left(1-\frac{1}{p}\right)\right)+\sum_{z\in Z,z\neq prime}\log \left(1+\frac{1}{z}\right)\\
&=\sum_p \log \left(1-\frac{1}{p^2}\right)+\sum_p \sum_{k=2}^\infty \log \left(1+\frac{1}{p+\cdots +p^k}\right)\\
&=\sum_p \log \left(1-\frac{1}{p^2}\right)+\sum_p \log\prod_{k=2}^\infty\left(\frac{1+p^{-1}+\cdots+p^{-k}}{1+p^{-1}+\cdots+p^{-(k-1)}}\right)\\
&=\sum_p \log \left(1-\frac{1}{p^2}\right)+\sum_p \log\prod_{k=2}^\infty\left(\frac{1-p^{-k-1}}{1-p^{-k}}\right)\\
&=\sum_p \left(\log \left(1-\frac{1}{p^2}\right)+\log \frac{1}{1-p^{-2}}\right)=0.\tag{3.2}
\end{align*}
\end{proof}

\noindent {\bfseries Theorem 4.}
\textit{
\begin{equation}\tag{4.1}
\sum_{i=1}^m\frac{1}{z_i}-W_1<\log \rho(n_m)<\sum_{i=1}^m\frac{1}{z_i}-W_1+\frac{1}{2z_m}.
\end{equation}
Hence
\begin{equation}\tag{4.2}
\log \rho(n_m)=\sum_{i=1}^m\frac{1}{z_i}-W_1+O\left(\frac{1}{z_m}\right).
\end{equation}
}
\begin{proof}
From construction (8), we have
\begin{equation}\tag{4.3}
\log \rho(n_m)=\sum_{i=1}^m\delta_i=\sum_{i=1}^m\log\left(1+\frac{1}{z_i}\right).
\end{equation}
Therefore
\begin{align*}
\log \rho(n_m)&-\sum_{i=1}^m\frac{1}{z_i}=\sum_{i=1}^m \left(\log\left(1+\frac{1}{z_i}\right)-\frac{1}{z_i}\right)\\
&=\sum_{i=1}^\infty\left(\log\left(1+\frac{1}{z_i}\right)-\frac{1}{z_i}\right)-\sum_{i=m+1}^\infty\left(\log\left(1+\frac{1}{z_i}\right)-\frac{1}{z_i}\right).\tag{4.4}
\end{align*}
The first term is $-W_1$. The second term has lower bound 0 and upper bound
\begin{align*}
-\sum_{i=m+11}^\infty\left(\log\left(1+\frac{1}{z_i}\right)-\frac{1}{z_i}\right)&=\sum_{i=m+1}^\infty \sum_{j=2}^\infty \frac{(-1)^j}{jz_i^j}\\
&<\sum_{i=m+1}^\infty \frac{1}{2z_i^2}<\frac{1}{2}\int_{z_m}^\infty \frac{1}{t^2}dt=\frac{1}{2z_m}.\tag{4.5}
\end{align*}
\end{proof}

Define
\begin{equation}
\psi_Z(x):=\sum_{z\in Z,z\leq x}\log z.
\end{equation}

\noindent {\bfseries Theorem 5.}
\textit{
\begin{multline}\tag{5.1}
\log \rho(n_m)=\log \log z_m+\gamma+\frac{\psi_Z(z_m)-z_m}{z_m \log z_m}\\
-\int_{z_m}^\infty \frac{(\psi_Z(y)-y)(\log y+1)}{y^2(\log y)^2}dy+O\left(\frac{1}{z_m}\right). 
\end{multline}
}
\begin{proof}
Using Stieltjes integral and integrating by parts, we have
\begin{align*}
\sum_{z\in Z,z\leq x}f(z)&=\int_{2-}^x \frac{f(y)}{\log y}d\psi_Z(y)\\
&=\frac{f(x)\psi_Z(x)}{\log x}-\int_2^x \psi_Z(y)\frac{d}{dy}\left(\frac{f(y)}{\log y}\right)dy\\
&=\frac{f(x)\psi_Z(x)}{\log x}-\int_2^xy\frac{d}{dy}\left(\frac{f(y)}{\log y}\right)dy+\int_2^xy\frac{d}{dy}\left(\frac{f(y)}{\log y}\right)dy\\
&\quad\quad-\int_2^x \psi_Z(y)\frac{d}{dy}\left(\frac{f(y)}{\log y}\right)dy\\
&=\frac{f(x)\psi_Z(x)}{\log x}-\frac{xf(x)}{\log x}+\frac{2f(2)}{\log 2}+\int_2^x\frac{f(y)}{\log y}dy\\
&\quad\quad-\int_2^x(\psi_Z(y)-y)\frac{d}{dy}\left(\frac{f(y)}{\log y}\right)dy\\
&=\frac{f(x)(\psi_Z(x)-x)}{\log x}+\frac{2f(2)}{\log 2}+\int_2^x\frac{f(y)}{\log y}dy\\
&\quad\quad-\int_2^\infty(\psi_Z(y)-y)\frac{d}{dy}\left(\frac{f(y)}{\log y}\right)dy\\
&\quad\quad+\int_x^\infty(\psi_Z(y)-y)\frac{d}{dy}\left(\frac{f(y)}{\log y}\right)dy.\tag{5.2}
\end{align*}
Substitute $f(x)=1/x$, we get
\begin{equation}\tag{5.3}
\frac{d}{dy}\left(\frac{f(y)}{\log y}\right)=\frac{d}{dy}\left(\frac{1}{y \log y}\right)=-\frac{\log y+1}{y^2(\log y)^2}.
\end{equation}
\begin{align*}
\sum_{z\in Z,z\leq x}\frac{1}{z}&=\frac{\psi_Z(x)-x}{x\log x}+\frac{1}{\log 2}+\int_2^x\frac{1}{x\log y}dy\\
&\quad +\int_2^\infty\frac{(\psi_Z(y)-y)(\log y+1)}{y^2(\log y)^2}dy-\int_x^\infty\frac{(\psi_Z(y)-y)(\log y+1)}{y^2(\log y)^2}dy\\
&=\frac{\psi_Z(x)-x}{x\log x}+\frac{1}{\log 2}+\log \log x-\log \log 2\\
&\quad +\int_2^\infty\frac{(\psi_Z(y)-y)(\log y+1)}{y^2(\log y)^2}dy-\int_x^\infty\frac{(\psi_Z(y)-y)(\log y+1)}{y^2(\log y)^2}dy\tag{5.4}
\end{align*}
We need to determine the constant. We have, by Mertens Theorem,
\begin{align*}
\frac{1}{\log 2}&-\log \log 2+\int_2^\infty\frac{(\psi_Z(y)-y)(\log y+1)}{y^2(\log y)^2}dy\\
&=\lim_{x\rightarrow\infty}\left(\sum_{z\in Z,z\leq x}\frac{1}{z}-\log \log x\right)\\
&=\lim_{x\rightarrow\infty}\left(\sum_{p\leq x}\frac{1}{p}-\log \log x+\sum_{z\in Z,z\leq x,x\neq prime}\frac{1}{z}\right)\\
&=M+\sum_{z\in Z,x\neq prime}\frac{1}{z}=W_2.\tag{5.5}
\end{align*}
Substitute (5.5)  into (5.4), we get
\begin{equation}\tag{5.6}
\sum_{z\in Z,z\leq x}\frac{1}{z}=\frac{\psi_Z(x)-x}{x\log x}+W_2+\log \log x-\int_x^\infty\frac{(\psi_Z(y)-y)(\log y+1)}{y^2(\log y)^2}dy.
\end{equation}
By theorems 3 and 4,
\begin{align*}
\log \rho (n_m)&<\sum_{i=1}^m\frac{1}{z_i}-W_1+O\left(\frac{1}{z_m}\right)\\
&<\frac{\psi_Z(z_m)-z_m}{z_m\log z_m}+W_2+\log \log z_m\\
&\quad-\int_{z_m}^\infty\frac{(\psi_Z(y)-y)(\log y+1)}{y^2(\log y)^2}dy-W_1+O\left(\frac{1}{z_m}\right)\\
&=\log \log z_m+\gamma+\frac{\psi_Z(z_m)-z_m}{z_m\log z_m}\\
&\quad-\int_{z_m}^\infty\frac{(\psi_Z(y)-y)(\log y+1)}{y^2(\log y)^2}dy+O\left(\frac{1}{z_m}\right).\tag{5.7}
\end{align*}
\end{proof}

 We now find bounds for $\psi_Z (x)$ in term of $\theta(x)$. 
 
\noindent {\bfseries Theorem 6.}
\textit{    Let $m\geq 2$ be an integer. Then for $z_m$ we have
\begin{equation}\tag{6.1}
\theta (z_m)+\sum_{k=2}^Kk\theta(z_m^{1/k})-\frac{(\log z_m)^2}{\log 2}\leq \psi_Z(z_m),
\end{equation}
\begin{equation}\tag{6.2}
\psi_Z(z_m)<\theta (z_m)+\sum_{k=2}^Kk\theta(z_m^{1/k})+\frac{2\log z_m \log \log z_m}{\log 2},
\end{equation}
where $K=\lfloor \log z_m /\log 2\rfloor$ (cf. Lemma 2).
\begin{equation}\tag{6.3}
\psi_Z(z_m)=\theta (z_m)+2\sqrt{z_m}+O\left(\frac{\sqrt{z_m}}{\log^2 z_m}\right).
\end{equation}
}
\begin{proof}
\begin{align*}
\psi_Z(z_m)&=\theta (z_m)+\sum_{k=2}^K\sum_{q\leq y_k}\log (q+\cdots+q^k)\\
&\geq \theta (z_m)+\sum_{k=2}^K\sum_{q\leq y_k}\log (q^k)\\
&=\theta (z_m)+\sum_{k=2}^K k\theta(y_k)\\
\text{(by Lemma 1)}&\geq \theta (z_m)+\sum_{k=2}^K k\left(\theta(z_m^{1/k})-\frac{\log z_m}{k}\right)\\
\text{(by Lemma 2)}&\geq \theta (z_m)+\sum_{k=2}^K k\theta(z_m^{1/k})-\left(\left\lfloor\frac{\log z_m}{\log 2}\right\rfloor-1\right)\log z_m\\
&\geq \theta (z_m)+\sum_{k=2}^K k\theta(z_m^{1/k})-\frac{(\log z_m)^2}{\log 2}.\tag{6.4}
\end{align*}
For (6.2), we have, by Mertens theorem,
\begin{align*}
\psi_Z(z_m)&=\theta (z_m)+\sum_{k=2}^K\sum_{q\leq y_k}\log (q+\cdots+q^k)\\
&=\theta (z_m)+\sum_{k=2}^K\sum_{q\leq y_k}\log (q^k(1+q^{-1}+\cdots+q^{-k+1}))\\
&<\theta (z_m)+\sum_{k=2}^K k\sum_{q\leq y_k}\log q+\sum_{k=2}^K\sum_{q\leq y_k}\log(1+2q^{-1})\\
&<\theta (z_m)+\sum_{k=2}^K k\theta(y_k)+2\sum_{k=2}^K\sum_{q\leq y_k}q^{-1}.\tag{6.5}
\end{align*}
Using Mertens Theorem for large $y_k$ and numerical calculation for small $y_k$, we have
\begin{equation}\tag{6.6}
\sum_{q\leq y_k}q^{-1} <\log \log y_k+0.8666.
\end{equation}
Hence
\begin{align*}
\psi_Z(z_m)&<\theta (z_m)+\sum_{k=2}^K k\theta(y_k)+2\sum_{k=2}^K(\log \log y_k+0.8666)\\
&<\theta (z_m)+\sum_{k=2}^K k\theta(y_k)+2\sum_{k=2}^K\log \log z_m +2\sum_{k=2}^K(-\log k+0.8666)\\
&<\theta (z_m)+\sum_{k=2}^K k\theta(z_m^{1/k})+\frac{2\log z_m \log \log z_m}{\log 2},\quad \forall \,m\geq 7.\tag{6.7}
\end{align*}
The cases of $2\leq m \leq 6$ can be numerically verified.

To prove (6.3), set $k=2, \eta_k=0.2, x_k=3\, 594\, 641$ in [Dusart 2018] Theorem 4.2, we have for $x\geq 3\, 594\, 641$,
\begin{equation}\tag{6.8}
\lvert \theta(x)-x\rvert<\frac{0.2x}{\log^2 x}.
\end{equation}
Hence
\begin{align*}
\psi_Z(z_m)&=\theta (z_m)+2\sqrt{z_m}+O\left(\frac{0.2\sqrt{z_m}}{\log^2 \sqrt{z_m}}\right)\\
&=\theta (z_m)+2\sqrt{z_m}+O\left(\frac{\sqrt{z_m}}{\log^2 z_m}\right).\tag{6.9}
\end{align*}
\end{proof}

\noindent {\bfseries Theorem 7.}
\textit{Let $m\geq 1$ be an integer. Then
\begin{equation}\tag{7.1}
\psi(z_m)-\log z_m\log \log z_m<\log n_m\leq \psi(z_m).
\end{equation}
}
\begin{proof}
Using notation in formulas (8) and (9), by Lemma 1,
\begin{equation}\tag{7.2}
\theta(z_m^{1/k})-\frac{1}{k}\log z_m\leq \theta(y_k)\leq \theta(z_m^{1/k}).
\end{equation}
For the right inequality
\begin{equation}\tag{7.3}
\log n_m=\sum_{k=1}^K\theta(y_k)<\sum_{k=1}^K \theta(z_m^{1/k})=\psi(z_m),
\end{equation}
where $K=\lfloor\log z_m/\log 2\rfloor$. For the left inequality, 
\begin{align*}
\log n_m&\geq\sum_{k=1}^K\theta(z_m^{1/k})-\sum_{k=2}^K \frac{1}{k}\log z_m\\
&=\psi(z_m)-\log z_m\sum_{k=2}^K\frac{1}{k}\\
&>\psi(z_m)-\log z_m(\log K+\gamma-1)\\
&\geq\psi(z_m)-\log z_m(\log \log z_m-\log 2+\gamma-1)\\
&>\psi(z_m)-\log z_m\log \log z_m.\tag{7.4}
\end{align*}
\end{proof}

\noindent {\bfseries Theorem 8.}
\textit{
\begin{align*}
\log \log \log n_m&=\log \log z_m+\frac{\psi(z_m)-z_m}{z_m \log z_m}\\
&\quad +O\left(\frac{(\psi(z_m)-z_m)^2}{z_m^2\log z_m}\right)+O\left(\frac{\log \log z_m}{z_m}\right).\tag{8.1}
\end{align*}
}
\begin{proof}
By Theorem 7,
\begin{align*}
\log &\log n_m-\log z_m\\
&=\log(\psi(z_m)+O(\log z_m\log \log z_m))-\log z_m\\
&=\log \psi(z_m)+O\left(\frac{\log z_m\log \log z_m}{\psi(z_m)}\right)-\log z_m\\
&=\log \left(\frac{\psi(z_m)}{z_m}\right)+O\left(\frac{\log z_m\log \log z_m}{z_m}\right)\\
&=\log \left(1+\frac{\psi(z_m)-z_m}{z_m}\right)+O\left(\frac{\log z_m\log \log z_m}{z_m}\right)\\
&=\frac{\psi(z_m)-z_m}{z_m}+O\left(\frac{(\psi(z_m)-z_m)^2}{z_m^2}\right)+O\left(\frac{(\log z_m\log \log z_m}{z_m}\right).\tag{8.2}
\end{align*}
Hence
\begin{align*}
\log &\log \log n_m-\log \log z_m=\log \left(1+\frac{\log \log n_m-\log z_m}{\log z_m}\right)\\
&=\sum_{i=1}^\infty \frac{(-1)^{i-1}}{i}\left(\frac{\log \log n_m-\log z_m}{\log z_m}\right)^i\\
&=\frac{\log \log n_m-\log z_m}{\log z_m}+O\left(\frac{(\log \log n_m-\log z_m)^2}{(\log z_m)^2}\right)\\
&=\frac{\psi(z_m)-z_m}{z_m\log z_m}+O\left(\frac{(\psi(z_m)-z_m)^2}{z_m^2\log z_m}\right)+O\left(\frac{\log \log z_m}{z_m}\right).\tag{8.3}
\end{align*}
\end{proof}

\noindent {\bfseries Lemma 3.}
\textit{
If the Riemann hypothesis is false, then there exists a real b with $0<b<1/2$ and a sequence $(x_i )$ of reals such that  $x_i$ is not a prime,
\begin{equation}\tag{L3.1}
x_i-\theta(x_i)=\frac{bx_i^{1-b}\log^2 x_i}{\log x_i+1} <bx_i^{1-b}\log x_i
\end{equation}
and
\begin{equation}\tag{L3.2}
\int_{x_i}^\infty \frac{(\theta(y)-y)(\log y+1)}{y^2\log^2 y}dy<-x_i^{-b}.
\end{equation}
}
\begin{proof}
Write
\begin{equation}\tag{L3.3}
K(x):=\int_x^\infty \frac{(\theta(y)-y)(\log y+1)}{y^2\log^2 y}dy.
\end{equation}
By Step (5) of the proof of Theorem 5.29 of [Broughan 2017], there exists a real b with $0<b<1/2$ such that
\begin{equation}\tag{L3.4}
K(x)=\Omega_\pm(x^{-b}),
\end{equation}
where $\Omega_\pm$ is the oscillation symbol. That means there exists a sequence $(x_i )$, $\lim_{i\rightarrow \infty}x_i=\infty,$ of reals such that
\begin{equation}\tag{L3.5}
K(x_i)\begin{cases}
<-x_i^{-b}, & \text{if i is odd}\\
>x_i^{-b}, & \text{if i is even}
\end{cases}.
\end{equation}
Fix an odd i. Then (L3.2) holds for $x_i$. We can choose $x_i$ so that $K(x_i )+x_i^{-b}$ is a local minimum of $K(x)+x^{-b}$. Define
\begin{equation}\tag{L3.6}
r(x):=\frac{x^2\log^2x}{\log x+1}\cdot \frac{d}{dx}(K(x)+x^{-b})=x-\theta(x)-\frac{bx^{1-b}\log^2x }{\log x+1}.
\end{equation}
Then r(x) changes sign at $x_i$. Assume $x_i$ lies in an interval $[p,p')$ of two consecutive primes. 

We claim that $x_i\neq p$. For if $x_i=p$, then
\begin{align*}
\lim_{x\rightarrow p-0}r(x)&=p-(\theta(p)-\log p)-\frac{bx^{1-b}\log^2x }{\log x+1}\\
&=\lim_{x\rightarrow p+0}r(x)+\log p>\lim_{x\rightarrow p+0}r(x).\tag{L3.7}
\end{align*}
If $r(x)$ changed sign at p, $r(x)$ would be positive on the left of $x_i$ and negative on the right. That is,  $K(p)+p^{-b}$ would be a local maximum, which contradicts $K(x_i )+x_i^{-b}$ being a local minimum.

    Therefore, we must have $x_i\in (p,p')$ and
\begin{equation}\tag{L3.8}
x_i-\theta(x_i)-\frac{bx_i^{1-b}\log^2x_i }{\log x_i+1}=0.
\end{equation}
Hence (L3.1) holds for all $x_i$ with odd i. The sequence $(x_i )$ in the lemma statement can be set to $(x_{2i+1} )$ above.
\end{proof}

\noindent {\bfseries Lemma 4.}
\textit{
If the Riemann hypothesis is false, then there exists a real b with $0<b<1/2$ and a sequence $(n_m )$ of LR numbers such that 
\begin{equation}\tag{L4.1}
\frac{bz_m^{1-b}\log^2 z_m}{\log z_m+1}-z_m^{0.525}\leq z_m-\theta(z_m)<bz_m^{1-b}\log z_m
\end{equation}
and
\begin{equation}\tag{L4.2}
\int_{z_m}^\infty \frac{(\theta(y)-y)(\log y+1)}{y^2\log^2 y}dy<-z_m^{-b}(1-z_m^{-0.475}).
\end{equation}
}
\begin{proof}
Write
\begin{equation}\tag{L4.3}
K(x):=\int_x^\infty \frac{(\theta(y)-y)(\log y+1)}{y^2\log^2 y}dy.
\end{equation}
By Lemma 3, there exists a sequence $(x_i )$ of reals such that $x_i$ is not a prime,
\begin{equation}\tag{L4.3}
x_i-\theta(x_i)=\frac{bx_i^{1-b}\log^2 x_i}{\log x_i+1}.
\end{equation}
and
\begin{equation}\tag{L4.4}
K(x_i)<-x_i^{-b}.
\end{equation}
Assume $x_i$ is in interval $[z_m,z_{m+1})\subset [p,p')$ with p and p' consecutive primes. If $x_i=z_m$, then we are done. So we may assume $x_i \in (z_m,z_{m+1})$. Define
\begin{equation}\tag{L4.5}
r(x):=\frac{x^2\log^2x}{\log x+1}\cdot \frac{d}{dx}(K(x)+x^{-b})=x-\theta(x)-\frac{bx^{1-b}\log^2x }{\log x+1}.
\end{equation}
\begin{equation}\tag{L4.6}
r'(x)=1-\frac{bx^{-b}\log x ((1-b)\log^2 x+(2-b)\log x+2}{(\log x+1)^2}.
\end{equation}
Then for $x^{-b} \log x<1$ we have
\begin{equation}\tag{L4.7}
0<r'(x)<1,
\end{equation}
By [BH 2001] Theorem 1, for sufficiently large x, there is at least one prime in $[x,x+x^{0.525} ]$. Hence 
\begin{equation}\tag{L4.8}
x_i-z_m<z_m^{0.525}.
\end{equation}
Since $r(x_i )=0$ and $r(x)<0$ for $x\in (z_m,x_i )$, we have
\begin{equation}\tag{L4.9}
0>r(z_m)>-\left(\max_{x\in (z_m,x_i)}r'(x)\right)z_m^{0.525}>-z_m^{0.525}.
\end{equation}
(L4.5) and (L4.9) mean (L4.1) is true. For (L4.2) we have
\begin{align*}
\left\lvert\int_{z_m}^{x_i} \frac{(\theta(y)-y)(\log y+1)}{y^2\log^2 y}dy \right\rvert &< (x_i-z_m)\frac{\lvert z_m-\theta(z_m)\rvert (\log z_m+1)}{z_m^2\log^2 z_m}\\
&<\frac{z_m^{0.525}(bz_m^{1-b}\log z_m+\log z_m)(\log z_m+1)}{z_m^2 \log^2 z_m}\\
&<z_m^{-0.475-b}.\tag{L4.10}
\end{align*}
Hence
\begin{align*}
\int_{z_m}^\infty &\frac{(\theta(y)-y)(\log y+1)}{y^2\log^2 y}dy\\
&= \int_{x_i}^\infty \frac{(\theta(y)-y)(\log y+1)}{y^2\log^2 y}dy+\int_{z_m}^{x_i} &\frac{(\theta(y)-y)(\log y+1)}{y^2\log^2 y}dy\\
&<-z_m^{-b}+z_m^{-0.475-b}=z_m^{-b}(1-z_m^{-0.475}).\tag{L4.11}
\end{align*}
So (L4.2) holds.
\end{proof}

\noindent {\bfseries Lemma 5.}
\textit{
\begin{equation}\tag{L5.1}
\int_x^\infty \frac{\sqrt{x}(\log y+1)}{y^2 \log^2 y}dy=O\left(\frac{1}{\sqrt{x}\log x}\right).
\end{equation}
}
\begin{proof}
By [Broughan 2017] Lemma 5.16 we have
\begin{align*}
\frac{2}{\sqrt{x}\log x}-\frac{2}{\sqrt{x}\log^2 x}&\leq \int_x^\infty \frac{\sqrt{x}(\log y+1)}{y^2 \log^2 y}dy\\
&\leq \frac{2}{\sqrt{x}\log x}-\frac{2}{\sqrt{x}\log^2 x}+\frac{8}{\sqrt{x}\log^3 x}.\tag{L5.2}
\end{align*}
So (L5.1) holds.
\end{proof}

\noindent {\bfseries Theorem 9.}
\textit{
Robin hypothesis is true if and only if all LR numbers $n_m>5040$ satisfy
\begin{equation}\tag{9.1}
\rho(n_m)<e^\gamma \log \log n_m.
\end{equation}
}
\begin{proof}
If Robin hypothesis is true, then (9.1) is true for all integers $n>5040$, hence true for all LR numbers $n_m>5040$.

    If Robin hypothesis is false, then Riemann hypothesis is false. By Theorem 6, we have
\begin{equation}\tag{9.2}
\psi_Z(y)=\theta(y)+2\sqrt{y}+O(\sqrt{y}/\log y),
\end{equation}
Hence by Theorems 5 and 8 and Lemma 5, we have, for large LR number $n_m$,
\begin{align*}
\log& \rho(n_m)-\log \log \log n_m-\gamma\\
&=-\int_{z_m}^\infty \frac{(\psi_Z(y)-y)(\log y+1)}{y^2 \log^2 y}dy\\
&\quad+O\left(\frac{(\psi(z_m)-z_m)^2}{z_m^2\log z_m}\right)+O\left(\frac{\log \log z_m}{z_m}\right)\\
&=-\int_{z_m}^\infty \frac{(\theta(y)-y)(\log y+1)}{y^2\log^2 y}dy\\
&\quad+O\left(\frac{(\theta(z_m)-z_m)^2}{z_m^2\log z_m}\right)+O\left(\frac{1}{\sqrt{z_m}\log z_m}\right).\tag{9.3}
\end{align*}
By Lemma 4, there exists a real b with $0<b<1/2$ and a sequence $(z_m )$ such that 
\begin{equation}\tag{9.4}
\frac{bz_m^{1-b}\log^2 z_m}{\log z_m+1}-z_m^{0.525}\leq z_m-\theta(z_m)<bz_m^{1-b}\log z_m
\end{equation}
and
\begin{equation}\tag{9.5}
-\int_{z_m}^\infty \frac{(\theta(y)-y)(1+\log y)}{y^2\log^2 y}dy>z_m^{-b}(1-z_m^{-0.475}).
\end{equation}
Substitute (9.4) and (9.5) into (9.3), we get
\begin{align*}
\log&\, \rho(n_m)-\log \log \log n_m-\gamma\\
&>z_m^{-b}+O\left(\frac{z_m^{1.05}}{z_m^2\log z_m}\right)+O\left(\frac{\log z_m}{z_m^{2b}}\right)+O\left(\frac{1}{\sqrt{z_m}\log z_m}\right)>0.\tag{9.6}
\end{align*}
That is, Robin inequality fails for some large $n_m$.
\end{proof}

\begin{center}
{\bfseries \large References}
\end{center}

\noindent {[}BH 2001{]} R. C. Baker and G. Harman, The difference between consecutive primes, II, Proc. London Math. Soc. (3) 83 (2001) 532-562.\\
{[}Broughan 2017{]} K. Broughan, \textit{Equivalents of the Riemann Hypothesis} Vol 1. Cambridge Univ. Press. (2017)\\
{[}CLMS 2007{]} Y. -J. Choie, N. Lichiardopol, P. Moree, and P. Sol\'{e}. \textit{On Robin’s criterion for the Riemann hypothesis}. J. Th\'{e}or. Nombres Bordeaux, 19(2):357–372, (2007).\\
{[}Dusart 2018{]} P. Dusart. \textit{Explicit estimates of some functions over primes}. Ramanujan J., 45(1):227–251, (2018).\\
{[}Robin 1984{]} G. Robin. \textit{Grandes valeurs de la fonction somme des diviseurs et hypoth\'{e}se de Riemann}. Journal de mathématiques pures et appliqu\'{e}es. (9), 63(2):187–213, (1984).

\end{document}